\documentclass[5p,times,numbers]{elsarticle}

\journal{TBA}

\usepackage{amsthm, amsmath, amssymb, amsfonts, graphicx, epsfig}
\usepackage{accents}

\usepackage{bbm}
\usepackage{mathrsfs,dsfont}
\usepackage{multirow}
\usepackage{booktabs,nicefrac}
\usepackage{array}
\usepackage{caption}
\usepackage[normalem]{ulem}
\usepackage{cancel}

\usepackage{epstopdf}

\usepackage{graphicx}
\usepackage[font=sl,labelfont=bf]{caption}
\usepackage{subcaption}

\usepackage{float}
\restylefloat{table}

\usepackage{enumerate}
\usepackage[colorlinks=true, pdfstartview=FitV, linkcolor=blue,
            citecolor=blue, urlcolor=blue]{hyperref}
\usepackage[usenames]{color} 

\usepackage[numbers]{natbib}

\usepackage{url}
\makeatletter\def\url@leostyle{%
 \@ifundefined{selectfont}{\def\UrlFont{\sf}}{\def\UrlFont{\scriptsize\ttfamily}}} \makeatother

\newtheorem{theorem}{Theorem}

\theoremstyle{definition}

\theoremstyle{remark}
\newtheorem{remark}[theorem]{Remark}

\numberwithin{equation}{section}
\numberwithin{theorem}{section}

\definecolor{Red}{rgb}{0.8,0,0.1}

\def\cF{\mathcal{F}}

\def\bE{\mathbb{E}}

\def\bN{\mathbb{N}}

\def\bP{\mathbb{P}}

\def\bR{\mathbb{R}}

\DeclareMathOperator*{\argmax}{arg\,max} 

\makeatletter
\def\namedlabel#1#2{\begingroup
    #2%
    \def\@currentlabel{#2}%
    \phantomsection\label{#1}\endgroup
}
\makeatother

\usepackage{bigfoot}
\usepackage{tcolorbox}
\usepackage{lipsum}
\usepackage{placeins}

\usepackage[british]{babel}

\usepackage{todonotes}
\usepackage{colortbl}

\begin{document}

\begin{frontmatter} 
\title{Blackwell optimality in risk-sensitive stochastic control$^{1,2}$}

\author[a1]{Marcin Pitera\fnref{fn1}}
\ead{marcin.pitera@uj.edu.pl}
\address[a1]{Institute of Mathematics, Jagiellonian University, Krakow, Poland}

\author[a3]{\L{}ukasz Stettner\fnref{fn1}}
\address[a3]{Institute of Mathematics, Polish Academy of Sciences, Warsaw, Poland}

\fntext[fn3]{This is the pre-submission version of a paper published in the proceedings of the 11th International Conference on Control, Decision and Information Technologies (CoDIT), 2025, Split, Croatia, $\copyright$ IEEE. The final version is available at https://doi.org/10.1109/CoDIT66093.2025.11321525.}

\fntext[fn2]{$\copyright$ 2026 IEEE. Personal use of this material is permitted.  Permission from IEEE must be obtained for all other uses, in any current or future media, including reprinting/republishing this material for advertising or promotional purposes, creating new collective works, for resale or redistribution to servers or lists, or reuse of any copyrighted component of this work in other works.}

\fntext[fn1]{Marcin Pitera and {\L}ukasz Stettner acknowledge research support by Polish National Science Centre grant no. 2024/53/B/HS4/00433. }

\begin{abstract}
In this paper, we consider a discrete-time Markov Decision Process (MDP) on a finite state-action space with a long-run risk-sensitive criterion used as the objective function. We discuss the concept of Blackwell optimality and comment on intricacies which arise when the risk-neutral expectation is replaced by the risk-sensitive entropy. Also, we show the relation between the Blackwell optimality and ultimate stationarity and provide an illustrative example that helps to better understand the structural difference between these two concepts.
\end{abstract}






\end{frontmatter}

\section{Introduction}

The {\em risk-sensitive criterion} in its {\it averaged} and {\it discounted} form is a popular discrete-time MDP objective function choice if one wants to incorporate risk into the decision-making process and consider long-run time horizon; see \cite{Whi1990,BauRie2011,BauJas2024} and references therein. Essentially, this criterion is a non-linear extension of the traditional risk-neutral mean criterion, in which the risk sensitivity parameter $\gamma\neq 0$ is encoded in the entropic-utility function $\frac{1}{\gamma}\ln \bE\left[e^{\gamma S} \right]$ that is applied to the cumulative reward $S$. Determining whether to use the averaged approach, where risk is considered per unit of time before eventually applying a limit, or the discounted criterion, which involves discounting future returns to ensure finite outcomes, is not straightforward and often relies on the underlying problem specification.

In this paper, we focus on the {\em Blackwell optimality}, a specific strategy property that allows transfer from the discounted to the averaged setup, see~\cite{Bla1962,Put2009}. On the one hand, this property holds significance in risk-neutral finite state-action framework, particularly within theoretical algorithms relying on the vanishing-discount approach and numerous Reinforced Learning techniques that use approximate discounting mechanisms, see~\cite{PorCavCru2023,BisBor2023,DewDunEshGalRoo2020,DewGal2022,GraPet2023} and references therein. On the other hand, it has been only recently shown in \cite{BauPitSte2024} that the Blackwell property can be efficiently extended from the risk-neutral setup to the risk-sensitive setup. The transfer is unexpectedly challenging due to the structural difference of the risk-neutral and risk-sensitive setups, even in the finite framework. In this paper, we investigate the sources of this challenge, recall recent results obtained in \cite{BauPitSte2024}, show how they are connected to the classical theory, and commenKalt how Blackwell optimality is linked to {\it ultimate stationarity} (also called {\it eventual-stationarity}), a concept introduced in \cite{Jaq1976} that links optimal policies in the risk-sensitive setup with optimal policies in the risk-neutral setup. For completeness, we also introduce an example in which the difference between those two concepts is illustrated. We hope that improved comprehension of Blackwell optimality within risk-sensitive framework may foster advancements in developing efficient risk-sensitive reinforcement learning algorithms, an area under active and intensive exploration, see \cite{FeiYanWan2021,DinJinLav2023,NooBar2021}, and references therein. 

\section{Long-run risk-sensitive stochastic control for MDPs}
Consider a finite state-action space $(E,U)$ and a time-homogeneous Markov Decision Process (MDP) given via the family of transition probabilities $(\bP^{a})_{a\in U}$.  As usual, we use $(\Omega,\cF,(\cF_n)_{n\in\bN})$ to denote the discrete-time filtered canonical probabability space and use $\Pi$ to denote the set of all policies, that is, sequences of random variables $\pi=(a_n)_{n\in\bN}$ such that $a_n\colon\Omega\to U$ is $\cF_n$-measurable. Furthermore, for $\pi\in\Pi$ and $x\in E$ we use $\bP_x^{\pi}$ to denote the controlled probability for the canonical process $(X_n)_{n\in\bN}$ with $X_0=x$. We also use $\Pi'$ and $\Pi''$ to denote the sets of all stationary Markov policies and  Markov policies, respectively. With a slight abuse of notation we link elements of $\Pi'$ with functions $u\colon E\to U$ and sometimes write $u\in \Pi'$, having in mind that $u$ constitutes policy $\pi=(a_n)_{n\in\bN}\in \Pi$ such that $a_n=u(X_n)$; we also write $\bP^u$ instead of $\bP^{\pi}$. Similarly, we associate elements of $\Pi''$ with sequences of functions $(u_n)_{n\in\bN}$. We refer to \cite{BauRie2011} and \cite{BauPitSte2024} for MDP construction details and formal definitions.

Given a fixed running reward function $c\colon E\times U\to\bR$, we consider the {\it long-run risk-sensitive averaged criterion}, that is, the objective function given by
\begin{equation}\label{eq:RSC}
J_\gamma(x,\pi):=
\begin{cases}
\liminf_{n\to\infty}\tfrac{1}{n}\frac{1}{\gamma }\ln\bE_x^{\pi}\left[e^{\gamma\sum_{i=0}^{n-1}c(X_i,a_i)}\right],& \textrm{if } \gamma\neq 0,\\
\liminf_{n\to\infty}\frac{1}{n}\bE_x^{\pi}\left[\sum_{i=0}^{n-1}c(X_i,a_i)\right],& \textrm{if }\gamma=0,
\end{cases}
\end{equation}
where $\bE_x^{\pi}$ denotes the expectation for $\bP_x^{\pi}$, $x\in E$, $\pi\in\Pi$, and $\gamma\in\bR$ is the risk-aversion parameter; for $\gamma=0$ we recover the classical risk-neutral criterion, while for $\gamma\neq 0$, we are in the risk-sensitive setup (see~\cite{Whi1990}). Our goal is to find optimal policy $\pi^*$ for which we have $J_{\gamma}(x,\pi^*)=\sup_{\pi\in \Pi}J_{\gamma}(x,\pi)$.

For simplicity, in this paper we decided to impose a {\it strong uniform mixing} assumption on the underlying MDP, that is, assume condition
\begin{equation}\tag{C}\label{C}
\inf_{a\in U}\inf_{x,y\in E}\bP^{a}(x,y)>0.   
\end{equation}
Assumption \eqref{C} implies the existence of a unique full-support invariant measure for any stationary Markov control and also ensures one-step transition equivalence for all controlled transition kernels. Although the results presented in this paper are true under much weaker assumptions linked to classical {\it mixing} and {\it multi-step communication}, we decided to present them in the simplified setup; see \cite{BauPitSte2024} details. Notably, under \eqref{C}, the optimal value of the objective function is independent of the starting point; see~\cite{BauPitSte2024} for details.

The value of the objective criterion defined in \eqref{eq:RSC} as well as the respective optimal strategies are often approximated using the discounted setup and so called vanishing discount approach, see \cite{CavHer2017}. For a given discount factor $\beta\in (0,1)$, the discounted analog of criterion \eqref{eq:RSC} is given by
\begin{equation}\label{eq:RSC:discounted}
J_\gamma(x,\pi;\beta):=
\begin{cases}
\frac{1}{\gamma }\ln\bE_x^{\pi}\left[e^{\gamma\sum_{i=0}^{n-1}\beta^ic(X_i,a_i)}\right],& \textrm{if } \gamma\neq 0,\\
\bE_x^{\pi}\left[\sum_{i=0}^{n-1}\beta^i c(X_i,a_i)\right],& \textrm{if }\gamma=0.
\end{cases}
\end{equation}
For consistency, we also set $J_\gamma(x,\pi;1)=J_\gamma(x,\pi)$. In contrast to \eqref{eq:RSC}, the optimal strategy for \eqref{eq:RSC:discounted} with $\gamma\neq 0$ does not need to be Markov stationary even under assumption \eqref{C}. In this paper, we study how the optimal strategies induced by \eqref{eq:RSC:discounted} are linked to the optimal strategies induced by \eqref{eq:RSC} and differentiate the concept of {\it Blackwell optimality} and {\it ultimate stationarity} that were introduced in \cite{Bla1962} and \cite{Jaq1973}, respectively. 

For brevity, we say that a policy is {\it average optimal} (or optimal for averaged problem), if it is optimal for \eqref{eq:RSC} and {\it $\beta$-discount optimal} (or optimal for $\beta$-discount problem), if it is optimal for \eqref{eq:RSC:discounted}.

\section{Blackwell optimality}
In a nutshell, a policy $\pi$ is said to be {\it Blackwell optimal}, if there exists $\beta_0<1$ such that $\pi$ is $\beta$-discount optimal for any $\beta\in (\beta_0,1)$; we also say that the policy $\pi$ satisfies the {\it Blackwell property}, if it is Blackwell optimal. In this section we recall the result from \cite{Bla1962} in which existence of Blackwell optimal policies was shown, and later show how to extend this concept to risk-sensitive setup based on the recent results presented in \cite{BauPitSte2024}.

 \subsection{Risk-neutral setup}\label{S:rn}

In this subsection we fix $\gamma=0$ and consider risk-neutral setup. First, we recall the classical Blackwell property theorem which ensures existence of Blackwell optimal strategy in the finite setup

\begin{theorem}\label{th:blackwell}
Assume \eqref{C} and fix $\gamma=0$.  Then, there exists a stationary Markov policy $\pi\in\Pi'$ that satisfies the Blackwell property. Furthermore, policy $\pi$ is optimal for the risk-neutral averaged problem. 
\end{theorem}

The first proof of Theorem~\ref{th:blackwell} was presented in \cite{Bla1962}, see also Section 10.1.2 in \cite{Put2009} for an extensive discussion about Blackwell optimality and related proof concepts. The traditional approach relies on Laurent series expansions along with Cramer's rule to demonstrate that every element of the inverse matrix $(I-\beta \bP^u)^{-1}$ is a rational function of $\beta$, where $\bP^u$ represents the transition matrix under decision rule $u$. This approach is then coupled with the insight that the optimal policy for any given $\beta$ can be identified among a fixed set of stationary policies, thereby reducing the problem to a finite set of policies. This also shows that the interval $[0,1)$ for the discount factor can be partitioned into a finite number of segments, each associated with a decision rule that specifies the optimal stationary policy for those discount factors.

The concept of Blackwell optimality is now considered as a canonical result in the risk-neutral criterion MDP literature, and was extensively studied in numerous papers. We refer to \cite{DekHor1988,hordijk2002blackwell,cavazos1999direct,Put2009} for details on early developments in this area and to \cite{DewGal2022,GraPetVie2025} for a more recent contributions, which include comprehensive discussion of the Blackwell property as well as its various modifications. Also, as noted in \cite{GraPet2023}, the concept of Blackwell optimality is an important topic in Reinforced Learning for average reward optimality, as it allows a transition from discounter to average setup. In fact, better understanding of the Blackwell optimality criterion is stated as {\it one of the pressing questions} in the list of open research problems formulated in \cite{DewDunEshGalRoo2020}.

 \subsection{Risk-sensitive setup}

Unfortunately, it is not possible to directly transfer Theorem~\ref{th:blackwell} and the previously discussed proof methods to the risk-sensitive setup. In contrast to the optimal discounted policies for risk-neutral setup, the risk-sensitive discounted optimal policies usually lack stationarity, resulting in a non-linear challenge that inhibits direct use of techniques based on Laurent series expansion.  However, a modified Blackwell property can still be retrieved within the non-stationary framework through proofs relying on a blend of vanishing discount and span-contraction methods.

Consider a constant risk-sensitivity parameter $\gamma\neq 0$ with the optimal $\beta$-discounted non-stationary Markov policy represented as $(u_0^\beta, u_1^\beta, u_2^\beta,\ldots)$. Then, for each $n\in\bN$, there exists $\beta_n<1$ such that for any $\beta$ within $(\beta_n,1)$ the stationary Markov policy $(u_n^\beta,u_n^\beta,\ldots )$ is optimal for the risk-sensitive average problem.

To make this statement more formal, let us consider the optimal policy induced by the discounted Bellman equation for problem \eqref{eq:RSC:discounted}. Namely, for any $\beta\in (0,1)$ and $\gamma\neq 0$, let $\hat\pi^{\beta}:=(\hat{u}^{\beta}_0,\hat{u}^{\beta}_1,\ldots)$, constitute the optimal $\beta$-discount policy given by
\begin{equation}\label{eq:discounted.stat}
 \hat u^{\beta}_n(x):=\argmax_{a\in U}\left[c(x,a) + \frac{1}{\gamma_n}\ln \left( \sum_{y\in E} e^{w^\beta(y,\gamma_{n+1})}\mathbb{P}^a(x,y)\right)\right],
\end{equation}
where $\gamma_n:=\gamma\beta^n$ and $w^{\beta}(\cdot,\cdot)$ is the solution to the discounted Bellman equation
\begin{equation}\label{eq:Bellman2}
w^{\beta}(x, \gamma)=\max_{a\in U}\left[ c(x,a)+\frac{1}{\gamma}\ln\sum_{y\in E}e^{w^{\beta}(y,\beta \gamma)}\bP^a(x,y)\right].
\end{equation}
The policy $\hat\pi^{\beta}$ might be non-stationary as we might have $w(\cdot,\gamma_{n+1})\not\equiv w(\cdot,\gamma_{m+1})$ and consequently $\hat u^{\beta}_n\not\equiv \hat u^{\beta}_m$, for $n\neq m$; note that the Bellman equation corresponding to \eqref{eq:Bellman2} is in fact recursive and requires us to compute $w^{\beta}$ on parameter grid $(\gamma_n)_{n\in\bN}$. Also, one can show that $w^{\beta}(x,\gamma)=\sup_{\pi\in\Pi}J_{\gamma}(x,\pi;\beta)$. We are now ready to present the analogue of Theorem~\ref{th:blackwell} transferred to the risk-sensitive framework. 

\begin{theorem}\label{th:blackwell.sensitive}
Assume \eqref{C} and fix $\gamma\neq 0$. Then, for any $n\in\bN$ there is $\beta_n\in (0,1)$ such that for any $\beta\in(\beta_n,1)$ the stationary Markov policy $\hat u^{\beta}_n\in \Pi'$ defined in \eqref{eq:discounted.stat} is optimal for the risk-sensitive averaged problem.
\end{theorem}

The proof of Theorem~\ref{th:blackwell.sensitive} can be found in \cite{BauPitSte2024}. The proof techniques used therein are essentially different from those used to prove Theorem \ref{th:blackwell} and rely on a combination of span-contraction methods with vanishing discount approach. One of the key results (also presented and proved in \cite{BauPitSte2024}) which contributes to the proof of Theorem~\ref{th:blackwell.sensitive} is the vanishing discount approximation result, which states that for appropriately chosen sequences we can recover the solution to the averaged Bellman equation from the discounted value functions. Namely, let us fix $z\in E$ and consider centered function $\bar w^{\beta}(x,\gamma) :=w^{\beta}(x,\gamma)-w^{\beta}(z,\gamma)$. Furthermore, for $n\in\bN$, let
\begin{align*}
\lambda_n^{\beta}(\gamma) &:= w^{\beta}(z,\gamma\beta^n)- w^{\beta}(z,\gamma\beta^{n+1}),\\
\bar w_n^{\beta}(x,\gamma) &:=\bar w^{\beta}(x,\gamma\beta^n).
\end{align*}
Then, we can formulate the following theorem.
\begin{theorem}\label{th:add}
Assume \eqref{C}; fix $\gamma\neq 0$. Then, for $x\in E$ and $\gamma\in \bR\setminus\{0\}$, there exists limits
\[
w(x,\gamma):=\lim_{\beta\uparrow1} \bar w_n^{\beta}(x,\gamma) \,\textrm{ and }\, \lambda(\gamma):=\lim_{\beta\uparrow 1}\lambda_{n}^{\beta}(\gamma).
\]
Furthermore, the function $w(\cdot,\gamma)$ and the constant $\lambda(\gamma)$ constitute the optimal stationary policy to the averaged risk-sensitive problem, that is, Markov stationary policy
\[
u(x):=\argmax_{a\in U}\left[c(x,a)-\lambda(\gamma)+\frac{1}{\gamma}\ln\sum_{y\in E}e^{\gamma w(y,\gamma)}\bP^a(x,y)\right],
\]
is optimal for the averaged risk-sensitive criterion~\eqref{eq:RSC}.
\end{theorem}
 Theorem~\ref{th:add} proof follows directly from the proof of Theorem 4.4 presented in \cite{BauPitSte2024}. The proof is based on recursive selection schemes and use finite space induced uniform bounds for value functions and optimal constants. We refer to Proposition 4.2, Proposition 4.3, and proof of Theorem 4.4 in \cite{BauPitSte2024} for more details.

As in the risk-neutral setup, the Bellman property could be used to provide a link between optimal discounted and optimal averaged policy for the risk-sensitive criterion. This property is important in the context of theoretical algorithms based on the vanishing-discount approach as well as many Machine Learning methods based on approximating discounting schemes, see~\cite{CavHer2011,PorCavCru2023,BisBor2023,DewDunEshGalRoo2020}  and references therein. Consequently, we believe that better understanding of Blackwell optimality in the risk-sensitive setup could help to boost development of efficient risk-sensitive reinforced learning algorithms which is an area of on-going and intensive research, see~\cite{FeiYanWan2021,DinJinLav2023,NooBar2021} and references therein.

\section{Structural difference between risk-neutral and risk-sensitive setup}

In this section, we explain in details the structural difference between Blackwell optimality  in the risk-neutral and risk-sensitive setup. They key difference is linked to the fact that while optimal discount policies in the risk-neutral setup induced by Bellman equation are inherently stationary, this is not the case for the risk-sensitive setup. Consequently, it is structurally impossible to directly transfer Theorem~\ref{th:blackwell} to the risk-sensitive framework, and one should adjust for non-stationary policies, as done in Theorem~\ref{th:blackwell.sensitive}.

To better understand this, let us recall that the value function in the risk-sensitive discounted Bellman equation given in \eqref{eq:Bellman2}, that is, the value of $w^{\beta}(x, \gamma)$ in equation
\begin{equation}\label{eq:Bellman2b}
w^{\beta}(x, \gamma)=\max_{a\in U}\left[ c(x,a)+\frac{1}{\gamma}\ln\sum_{y\in E}e^{w^{\beta}(y,\beta \gamma)}\bP^a(x,y)\right]
\end{equation}
might strongly depend on the underlying risk sensitivity parameter $\gamma$ and consequently we might get $\hat u^{\beta}_n\not\equiv \hat u^{\beta}_m$, for $n\neq m$, where $\hat u^{\beta}_n$ is defined in \eqref{eq:discounted.stat}. On the other hand, when $\gamma=0$, the Bellman equation \eqref{eq:Bellman2b} simplifies to equation
\[
w^{\beta}(x, 0)=\max_{a\in U}\left[ c(x,a)+\beta\sum_{y\in E}w^{\beta}(y,0)\bP^a(x,y)\right],
\]
so that we only need to consider the value $w^{\beta}(\cdot, 0)$ and the respective stationary Markov policy $\hat u^{\beta}\in\Pi'$ given by
\begin{equation}\label{eq:discountedn.stat} 
 \hat u^{\beta}(x):=\argmax_{a\in U}\left[c(x,a) + \beta\left( \sum_{y\in E} w^\beta(y,0)\mathbb{P}^a(x,y)\right)\right].
\end{equation}
Essentially, this is due to the fact that the expectation operator embedded in the sum in \eqref{eq:Bellman2b} is positively homogeneous, while the entropy operator encoded in the logarithm of the expectation of the exponent in \ref{eq:Bellman2b} is a proper convex (or concave) function which cannot be linearly scaled. Furthermore, it should be noted that while the set of stationary Markov policies $\Pi'$ is finite, the set of Markov policies $\Pi''$ is not, which inhibits the usage of methods described in Section~\ref{S:rn}. This is essentially why the equivalent of Theorem~\ref{th:blackwell} cannot be true in the risk-sensitive setup and why the related proof techniques do not work in this setup. Namely, we cannot restrict our analysis to the finite number of stationary Markov policies in the risk-sensitive cane. On the other hand, this raises a natural question under which additional conditions there exists a stationary solution to the risk-sensitive discounted problem for any discount factor. In particular, one expects this to be the case for values of $\gamma$ which are sufficiently close to zero. In this case, the entropy operator becomes almost linear (as, in the limit, it converges to the expectation operator). This simple idea leads to the concept of {\it ultimate stationarity} which is formulated in the next subsection.

\section{Interaction between blackwell optimality and ultimate stationarity in the risk-sensitive setup}

The concept of {\it ultimate stationarity} was initially developed in \cite{Jaq1976}, where a specific form of risk-neutral optimal strategy, also optimal in the risk-aversion neighborhood of zero, was studied. We say that $\pi\in \Pi''$ is an {\it ultimately  stationary policy}, if it can be represented as
\begin{equation}\label{eq:ult.stat}
\pi=(u_0,u_1,\ldots,u_{N-1},u,u,\ldots),
\end{equation} for some $N\in\bN$ and one-step policies $u_1,\ldots,u_{N-1},u\in\Pi'$. In other words, the ultimate stationary policy $\pi$ becomes a stationary Markov policy after a fixed number of steps. With a slight abuse of notation, we also refer to policy $u$ in \eqref{eq:ult.stat} as the ultimate stationary policy. The first (essential) step in showing the existence of an ultimate stationarity policy is formulated in the following theorem.

\begin{theorem}\label{th:jacquette1}
Assume \eqref{C} and fix $\beta\in (0,1)$. Then, there is $\gamma_{0}<0$ (resp. $\gamma_0>0$) and $u\in\Pi'$ such that $u$ is optimal for the discounted problem \eqref{eq:RSC:discounted} for any $\gamma\in [\gamma_0,0]$ (resp. $\gamma\in [0,\gamma_0]$).
\end{theorem}

The proof of Theorem~\ref{th:jacquette1} can be found in~\cite{Jaq1973}. Knowing this result, one can deduce that, for fixed $\beta\in (0,1)$ and $\gamma\in\bR\setminus\{0\}$, the $n$th one-step policy in \eqref{eq:discounted.stat} could be seen as the $1$st one-step policy for the discounted problem with initial risk-aversion equal to $\beta^n\gamma$, which is closer to zero than the initial parameter $\gamma$. Note that ultimately stationary policy $u$ in representation~\eqref{eq:ult.stat} corresponds to policy from Theorem~\ref{th:jacquette1}. This observation can be used to prove the existence of an ultimate stationary policy for any value of $\beta\in (0,1)$; of course, one expects that the closer the value of $\beta$ to 1 is, the more steps one needs to make to ensure the stationarity property in the optimal policy. For details, we refer to \cite{Jaq1976} where this result is formulated and proved.

Upon initial inspection, one would hope that we can use Theorem~\ref{th:jacquette1} to somehow recover the statement of Theorem~\ref{th:blackwell} or Theorem~\ref{th:blackwell.sensitive}. Unfortunately, this is not the case, as those theorems are of different nature. Namely
\begin{itemize}
\item {\it The Blackwell policies} are linked to the situation in which risk aversion $\gamma$ is fixed and we let $\beta\to 1$;
\item {\it The ultimate stationary policies} are linked to the situation in which $\beta$ is fixed and we let  $\gamma\to 0$.
\end{itemize}
Those properties are structurally different and there is no direct connection between them. In particular, apart for some special situations, in which the initial risk-aversion parameter is very close to zero, the sets of Blackwell policies and ultimate stationary could be disjoint. This is presented in the illustrative example in Section~\ref{S:example}.

\begin{remark}
The existence of ultimate stationary policies (and Blackwell policies) is inherently linked to the finite space structure, in both risk-neutral and risk-sensitive cases. In fact, in the denumerable case, one can show direct existence counterexamples if no additional ergodic-mixing assumptions are imposed. We refer to \cite{Jaq1976} and \cite{BraCavFer1998} for details.
\end{remark}

\section{Illustrative example}\label{S:example}

In this section, we present an extended and refined analysis of the example introduced in~\cite{BauPitSte2024} that was based on~\cite{Jaq1976}. We use it to illustrate the structural difference between risk-neutral and risk-sensitive setups, and to visualise different nature of the Blackwell optimality and the ultimate stationarity property.

\begin{figure}[thpb]
      \centering
      \framebox{\parbox{3in}{\includegraphics[scale=0.5]{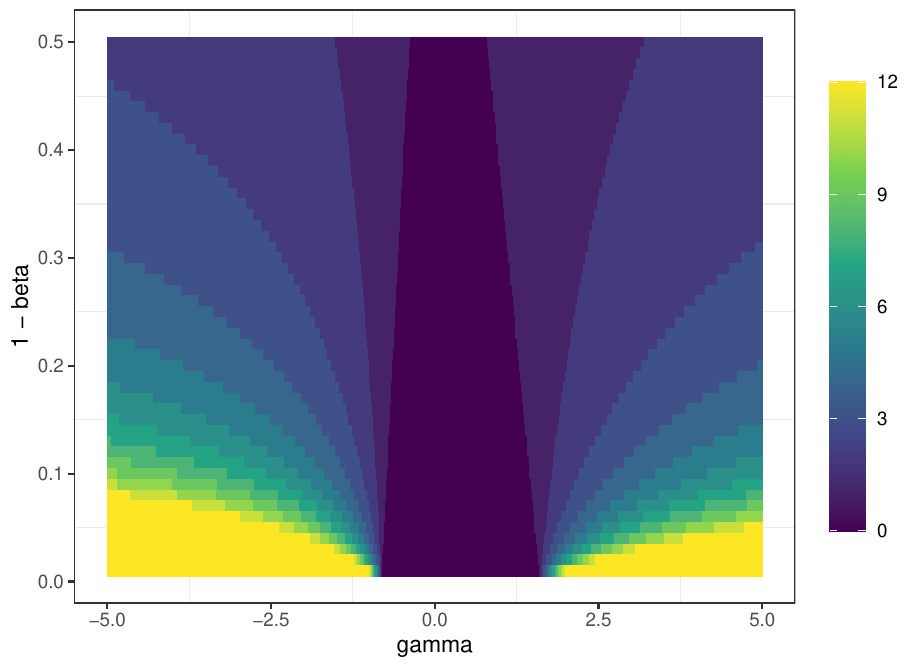}}}
      \caption{The plot presents the optimal policy switching point $n(\beta,\gamma)$ from Markov policy $u_0$ to Markov policy $u_1$ for different values of $\gamma\in [-5,5]$ and $\beta\in (0.5,1)$.}
      \label{F:ex1a}
   \end{figure}

Let $(E,U):=(\{1,2,3\},\{0,1\})$ with transition matricesgiven by
\[
P^0=\begin{bmatrix}
2\epsilon & 0.5-\epsilon  & 0.5-\epsilon \\
1-2\epsilon & \epsilon  & \epsilon \\
1-2\epsilon & \epsilon  & \epsilon \\
\end{bmatrix},
\]
\[
P^1=\begin{bmatrix}
2\epsilon & 0.9-\epsilon  & 0.1-\epsilon \\
1-2\epsilon & \epsilon & \epsilon \\
1-2\epsilon & \epsilon  & \epsilon \\
\end{bmatrix},
\]
for $\epsilon\in (0,0.1)$, and cost function $c\colon E\times U\to\bR$ defined as
\[
c(x,0)=
\begin{cases}
0& \textrm{if } x=1\\
0& \textrm{if } x=2\\
8& \textrm{if } x=3
\end{cases}\quad\textrm{and}\quad 
c(x,1)=
\begin{cases}
1& \textrm{if } x=1\\
0& \textrm{if } x=2\\
8& \textrm{if } x=3
\end{cases}.
\]
In this example, we have eight Markov stationary policies. Noting that $\bP^1(i,\cdot)\equiv \bP^2(i,\cdot)$ and  $c(i,\cdot)\equiv \textrm{const}$ for $i\in \{2,3\}$, we conclude that it does not matter which action is assigned to state $2$ and $3$. Consequently, without loss of generality, it is sufficient to consider Markov policies $\pi=(a_n)_{n\in\bN}$, where $a_n=u_0(X_n)$ or $a_n=u_1(X_n)$ for every $n\in\bN$, where $u_0(\cdot)\equiv0$ and $u_1(\cdot)\equiv 1$. Furthermore, for $\epsilon$ close to zero, following similar reasoning to the one introduced in \cite{BauPitSte2024}, one can conclude that for any $\beta\in (0,1)$ and $\gamma\in\bR$ the optimal Markov policy $\hat \pi(\beta,\gamma)\in\Pi''$ must be of the form
\[
\pi(\beta,\gamma)= (u_0,\ldots, u_0,u_1,u_1,\ldots),
\]
that is, in the fixed number of initial steps, say $n(\beta,\gamma)$, we use policy $u_0$, and then switch to policy $u_1$. As shown in \cite{BauPitSte2024}, the policy $u_0$ is typically Blackwell optimal and $u_1$ is ultimately stationary. To illustrate this, we have computed the switch value points, i.e. values of $n(\beta,\gamma)$ on the $(0.5,1)\times [-5,5]$ dyadic grid of step-size $\delta:=0.01$. For simplicity we performed the calculation for the limit case $\epsilon=0$ as it allows decomposing $(X_n)_{n\in\bN}$ into a series of independent (two-step) innovations which in turn allows decomposing entropic utility of the sum into entropic utilities of individual components. For brevity, we skip the technical calculation and focus on the output results. The obtained structure of optimal policies, represented by the value of $n(\beta,\gamma)$, is presented via a heat map in Figure~\ref{F:ex1a}.

The plot nicely illustrates the structural difference between {\it ultimate stationarity} and {\it Blackwell property}:

\begin{itemize}
    \item The {\it ultimate stationarity property} is represented by horizontal lines in Figure~\ref{F:ex1a}. Namely, for any fixed value of $\beta$, we ultimately switch to $u_1$ after $n(\beta,\gamma)$ steps. Furthermore, the closer we are to the risk-neutral case, the faster we switch; as expected, the policy $u_1$ becomes optimal stationary for values of $\gamma$ sufficiently close to zero.
    
   \item  The {\it Blackwell property} is represented by vertical lines in Figure~\ref{F:ex1a}. Namely, for any fixed value of $\gamma$ the number of steps $n(\beta,\gamma)$ goes either to $0$ or $\infty$ when $(1-\beta)\to 0$. The limit is determined by the optimal policy for the risk-sensitive averaged case: if $u_1$ is optimal for a given $\gamma$, then we get $n(\beta,\gamma)\to 0$ as $\beta\to 1$, and if $u_0$ is optimal, then we get $n(\beta,\gamma)\to \infty$ as $\beta\to 1$.
   \end{itemize}

\begin{figure}[thpb]
      \centering
      \framebox{\parbox{3in}{\includegraphics[scale=0.5]{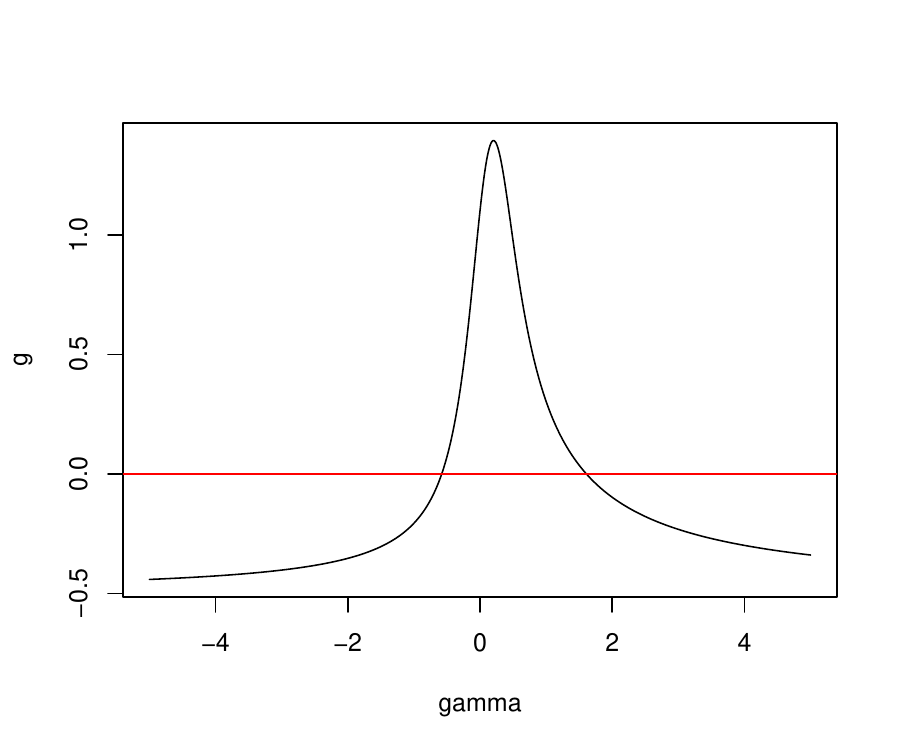}}}
      \caption{The plot presents the difference between $J_{\gamma}(1,u_1;1)$ and $J_{\gamma}(1,u_0;1)$  for different values of $\gamma\in [-5,5]$. Positive difference indicates optimality $u_1$, while negative difference indicates optimality of $u_0$.}
      \label{F:ex1b}
   \end{figure}

We note that the Blackwell optimal policy might be different than the ultimate stationary policy. In our example, this is the case for values of $\gamma$ sufficiently far away from zero. To better understand why this is the case, we can plot the difference between averaged risk-sensitive objective functions for $u_1$ and $u_2$, for different values of $\gamma\in\bR$. The plot of the function $g(\gamma):=[J_{\gamma}(1,u_1;1)-J_{\gamma}(1,u_0;1)]$ is presented in Figure~\ref{F:ex1b}.

 In a nutshell, Blackwell optimal policy recovers the risk averaged optimal policy which depends on the value of $\gamma\in\bR$. As expected, the switching points in Figure~\ref{F:ex1a} are matching the zeroes of function $g$ presented in Figure~\ref{F:ex1b}. This illustrates that while ultimate stationary policy is essentially recovering optimal behavior under almost risk-neutral risk regime, the Blackwell policy is recovering optimal behavior for the risk-sensitive case, when the discount is vanishing.

 {\footnotesize
\bibliographystyle{elsarticle-num}
 \bibliography{RSC_bibliografia}
 }

 \end{document}